\documentclass{amsart}
\usepackage{epsf,amssymb}

\begin{document}
\newcommand{\fx}{f^{(x)}}
\newcommand{\gx}{g^{(x)}}
\newcommand{\hx}{h^{(x)}}

\newcommand{\bP}{\bold P}
\newcommand{\bQ}{\bold Q}
\newcommand{\lessdoteq}{\preceq}
\renewcommand{\lessdot}{\prec}
\newcommand{\zeroh}{\hat 0}
\newcommand{\oneh}{\hat 1}
\newcommand{\bx}{\bf x}
\newcommand{\by}{\bf y}

\newtheorem{theorem}{Theorem}
\newtheorem{proposition}{Proposition}
\newtheorem{lemma}{Lemma}

\title [Graded left modular lattices]{Graded left modular lattices are supersolvable}
\author {Hugh Thomas}


\address {Fields Institute, 222 College Street, Toronto ON, M5T 3J1, Canada}
\curraddr{Department of Mathematics and Statistics, University
of New Brunswick, Fredericton NB, E3B 5A3 Canada}

\email {hugh@math.unb.ca}
\begin{abstract}
We provide a direct proof that a finite 
graded lattice with a 
maximal chain of left modular elements is supersolvable.  This result was
first 
established via a detour through EL-labellings in [MT] by combining results of 
McNamara [Mc] and Liu [Li].  As part of our proof, 
we show that the maximum graded quotient of the free product of a chain and
a single-element lattice is finite and distributive.
\end{abstract}

\maketitle

Supersolvability for lattices was introduced by Stanley [St].  A finite lattice
is {\it supersolvable} iff it has a maximal chain (called the $M$-chain) 
such that the sublattice generated
by the $M$-chain and any other chain is distributive.  

We say an element $x$ of a lattice is 
{\it left modular} if it satisfies:
$$(y\vee x)\wedge z=y\vee (x \wedge z)$$
for all $y\leq z$.   
Following Blass and Sagan [BS], 
we say that a lattice is left modular if it has a maximal chain of left
modular elements.  Stanley [St] showed that the elements of the $M$-chain of a 
supersolvable lattice are left modular, and thus that supersolvable lattices
are left modular.  

We say that a lattice is {\it graded} if, whenever $x<y$ and there is a 
finite maximal chain between $x$ and $y$, all the maximal chains
between $x$ and $y$ have the same length.  It is easy to check that 
supersolvable lattices are graded.  

The main result of our paper is the converse of these two results:

\begin{theorem} If $L$ is a finite, graded, left modular lattice,
then $L$ is supersolvable.  \end{theorem}

This result 
was first 
proved in [MT], as an immediate consequence of 
results of
Liu and McNamara.  
Liu [Li] showed that if a finite lattice is graded of rank $n$ 
and left modular, then it has an
EL-labelling of the edges of its Hasse diagram, 
such
that the labels which appear on any maximal chain are the numbers 1 through
 $n$ in some order.  McNamara [Mc] showed that for graded lattices of 
rank $n$, having such a labelling is equivalent to being supersolvable.  
These two results together immediately yield that finite graded left modular 
lattices are supersolvable.  However, since this proof involves considerations
which seem to be extraneous to the character of the 
result, it seemed worth giving a more direct and purely lattice-theoretic 
proof.  

On the way to our main result, we 
introduce the notion of the maximum graded quotient of a lattice.  
The maximum graded
quotient need not exist, but if it exists, it is unique.  We calculate 
explicitly the maximum graded quotient of the free product 
of the 
$k+1$-element chain $C_k$ with the single element lattice $S$ and show that
it is finite and distributive.

\section*{ The Maximum Graded Quotient of a Lattice}

When we refer to a quotient of a lattice, we mean a quotient with
respect to a lattice congruence, that is to say, a homomorphic image of the
original lattice.

Let $L$ be a lattice.  Define an equivalence relation $\sim$ on
$L$ by setting $x\sim y$ iff $\theta(x)=\theta(y)$ for all lattice 
homomorphisms $\theta$ from $L$ to a graded lattice.  It is straightforward
to check that $\sim$ is a lattice congruence.  We then define
$g(L)=L/\!\sim$.  By construction, $g(L)$ is the maximum quotient through
which every lattice homomorphism to a graded lattice factors.  

If $g(L)$ is graded, then we call it the maximum graded quotient of $L$.
Otherwise, we say that $L$ has no maximum graded quotient.  The lattice shown in 
Figure 1 has $g(L)=L$, and since $g(L)$ is not graded, $L$ has no
maximum graded quotient.

$$
\setlength{\unitlength}{3947sp}%
\begin{picture}(624,1074)(4789,-3598)
\thinlines
{\put(5101,-3511){\line(-1, 1){225}}
\put(4876,-3286){\line( 0, 1){450}}
\put(4876,-2836){\line( 1, 1){225}}
\put(5101,-2611){\line( 3,-4){225}}
\put(5326,-2911){\line( 0,-1){300}}
\put(5326,-3211){\line(-3,-4){225}}
}%
{\put(5101,-3511){\circle*{74}}
\put(4876,-3286){\circle*{74}}
\put(4876,-3061){\circle*{74}}
\put(4876,-2836){\circle*{74}}
\put(5101,-2611){\circle*{74}}
\put(5326,-2911){\circle*{74}}
\put(5326,-3211){\circle*{74}}}
\end{picture}
$$
$$\text{Figure 1}$$

For $x\in L$, we will write $[x]$ for the class of $x$ in $g(L)$.  
We write $a \lessdoteq b$ to indicate that either $a\lessdot b$ or 
$a=b$.

\begin{lemma} If $[x]\lessdoteq [y] \lessdoteq [z]$ (for instance,
if $x\lessdot y \lessdot z$ in $L$), and $[x] \leq [u] \leq [v] \leq [z]$, such
that $[u] \vee [y]= [z]$ and $[v] \wedge [y]= [x]$, then $[u]=[v]$.   
\end{lemma}

\begin{proof}  
 We consider separately graded quotients of $g(L)$ 
where $[y]$ is identified with $[x]$, where
$[y]$ is identified with $[z]$, where $[y]$ is not identified with either
$[x]$ or $[z]$, and where $[x]$, $[y]$, and $[z]$ are all identified.  
We see that 
in all these cases, $[u]$ and $[v]$ must be identified in the quotient.  
Since every lattice homomorphism from $L$ 
to a graded lattice factors through $g(L)$, this implies
that $u$ and $v$ are identified in any graded quotient of $L$, and therefore 
$[u]=[v]$.    
\end{proof}

\section*{ The Maximum Graded Quotient of $C_k \ast S$}

Let $C_k$ denote the chain of length $k$, with elements 
$x_0\lessdot\dots\lessdot x_k$.
Let $S$ denote the one element lattice, with a single element $y$.

\begin{lemma} The free product 
$C_k \ast S$ is a disjoint union of elements lying above
$x_0$ and elements lying below $y$.  \end{lemma}

\begin{proof} This is an immediate application of 
the Splitting Theorem [Gr, Theorem VI.1.11], which
says that the free product of two lattices $A$ and $B$ is the 
disjoint union of the dual ideal generated by $A$ and the ideal generated by
$B$.  
\end{proof}

We shall now proceed to consider these two subsets of $C_k\ast S$ in more
detail.  

\begin{lemma} The elements of $C_k \ast S$ lying below $y$
are exactly $y$ and $y\wedge x_i$ for $0\leq i\leq k$.  \end{lemma}

\begin{proof}
For $f\in C_k\ast S$, 
write $\fx$ for the smallest element of the $C_k$ which lies
above $f$.  If there is no such element, set $\fx=\oneh$.  We now claim that
$f\wedge y=\fx \wedge y$.

By definition, $\fx \geq f$, so $\fx\wedge y \geq f\wedge y$.  We prove the
other inequality by induction on the rank of a polynomial expression 
for $f$.  The statement is clearly true for rank 1 polynomials.  If the
rank of $f$ is greater than 1, it can be written as either $g\wedge h$
or $g \vee h$, for $g$ and $h$ polynomials of lower rank.  
Suppose
that $f=g\wedge h$.  Then $\fx=\gx \wedge \hx$ [Gr, Theorem VI.1.10], so
$$\fx \wedge y=\gx\wedge \hx \wedge y\leq g\wedge h \wedge y =f\wedge y.$$
Alternatively, 
suppose that $f=g\vee h$.  Then $\fx = \gx \vee \hx$ [Gr, Theorem VI.1.10]. 
Since 
$C_k \cup \{\oneh\}$ forms a chain, we may assume without loss of generality
that $\fx=\gx$.  Thus,
$$\fx \wedge y = \gx \wedge y \leq g \wedge y \leq (g\vee h)\wedge y = 
f\wedge y.$$
This completes the proof of the claim.  

It follows that if $z\leq y$, then 
$z=z\wedge y=z^{(x)}\wedge y$, and we have written
$z$ in the form described in the statement of Lemma 3. 
\end{proof}

\begin{lemma} The elements of $g(C_k \ast S)$ which lie strictly 
above $x_0$
are generated by $x_1,\dots$, $x_n$, $y\vee x_0$.  \end{lemma}

\begin{proof} We begin by showing that 
the elements of $C_k\ast S$ lying strictly above $x_0$ are generated by 
$x_1,\dots,x_n,y\vee x_0,(y\wedge x_1)\vee x_0,\dots,(y\wedge x_n)\vee x_0$.

Let $T_0$ denote $\{x_0,\dots,x_n,y\}$.  Define $T_i$ inductively as
those elements of $C_k \ast S$ which can be formed as either a meet or a 
join of a pair of 
elements in $T_{i-1}$.            The union of the $T_i$ is 
$C_k\ast S$.  We wish to show by induction on $i$ that any element of  
$T_i$ lying strictly above $x_0$
can be written as a polynomial in 
$x_1,\dots,x_n,y\vee x_0,(y\wedge x_1)\vee x_0,\dots,(y\wedge x_n)\vee x_0$.
The statement is certainly true for $i=0$.  Suppose it is true for $i-1$.
The statement is also true for any element of $T_i$ formed by a meet,
since if the meet lies strictly above $x_0$, so did both the elements of 
$T_{i-1}$.  Now consider the case of the join of two elements, $a$ and 
$b$, from $T_{i-1}$.  
If both $a$ and $b$ lie strictly above $x_0$, the statement is true for
 $a\vee b$ by induction.  
If neither $a$ nor $b$ lies strictly above $x_0$, then (by Lemma 3)
one of $a$ or $b$ must equal $x_0$, and $a\vee b$ is one of the generators
which we are allowing.  
Now suppose that $a$ lies strictly above $x_0$ and $b$ does not.  
By Lemma 3, $b$ equals
$x_0$, $y$, or $y\wedge x_i$.  If $b=x_0$, then $a\vee b=a$, and the statement
is true by induction.  Otherwise,  
$a\vee b=a\vee (b\vee x_0)$, and
$b\vee x_0$ is one of the allowed generators, so we are done.  
We have
shown that every element of $T_i$ lying above $x_0$ can be written in the 
desired form, and hence by induction that the same is true of any element 
of $C_k\ast S$ lying above $x_0$.

We now wish to show that the generators of the form 
$(y\wedge x_i)\vee x_0$ are unnecessary once we pass to $g(C_k\ast S)$.  
It follows from Lemma 3 that
$y\wedge x_n\lessdot y$.  Dually, $y\lessdot y\vee x_0$.  Observe that
$y\wedge x_n <(y\wedge x_n)\vee x_0 < (y\vee x_0)\wedge x_n < y\vee x_0$
in $C_k\ast S$.  Thus, by Lemma 1, 
$[(y\wedge x_n)\vee x_0]=[ (y\vee x_0)\wedge x_n]$.

We now proceed to show that 
$$[(y\wedge x_i)\vee x_0]=
[(y \vee x_0)\wedge x_i]$$ for all $1\leq i \leq n$.  
The proof is by downward induction;
we have already finished the base case, when $i=n$.  So suppose the result
holds for $i+1$.  In $L$,
$$y\wedge x_i \lessdot y\wedge x_{i+1} \lessdot (y\wedge x_{i+1})\vee x_0
<  (y\vee x_0)\wedge x_{i+1},$$ but when we pass to $g(L)$ the final
inequality becomes an equality by the induction hypothesis.  
Since in $L$ we also have 
that $$y\wedge x_i <(y \wedge x_i) \vee x_0<(y\vee x_0)\wedge x_i<
(y\vee x_0)\wedge x_{i+1},$$ we can apply Lemma 1 to conclude that 
$[(y\wedge x_i)\vee x_0]=
[(y \vee x_0)\wedge x_i]$ as desired.

We have already shown that the elements of $L$ 
lying above $x_0$ are generated by the $x_i$, $y\vee x_0$,
and the $(y\wedge x_i) \vee x_0$, for $i\geq 1$.  It follows that 
the elements of $g(L)$ above $[x_0]$ are
generated by
the $[x_i]$, $[y\vee x_0]$, and the $[( y\wedge x_i)\vee x_0]$.  But
$[(y\wedge x_i)\vee x_0]=[(y\vee x_0)\wedge x_i]=[y\vee x_0]\wedge [x_i]$,
and so the $[(y\wedge x_i)\vee x_0]$ are redundant, proving the lemma.
\end{proof}

\begin{proposition} The lattice $g(C_k \ast S)$ is as shown in 
Figure 2.  \end{proposition}

\setlength{\unitlength}{3947sp}%
\begingroup\makeatletter\ifx\SetFigFont\undefined%
\gdef\SetFigFont#1#2#3#4#5{%
  \reset@font\fontsize{#1}{#2pt}%
  \fontfamily{#3}\fontseries{#4}\fontshape{#5}%
  \selectfont}%
\fi\endgroup%
$$
\begin{picture}(2250,3066)(3526,-4894)
\put(3501,-4571){$x_0$}
\put(3501,-3971){$x_1$}
\put(3501,-2771){$x_{n-1}$}%
\put(3501,-2171){$x_n$}%
\thinlines
{\put(4201,-1861){\line(-1,-1){300}}
\put(3901,-2161){\line( 1,-1){300}}
\put(4201,-2461){\line(-1,-1){300}}
\put(3901,-2761){\line( 1,-1){300}}
}%
{\multiput(4201,-3061)(0.00000,-60){11}{\makebox(1.6667,11.6667){\SetFigFont{5}{6}{\rmdefault}{\mddefault}{\updefault}.}}
}%
{\multiput(4801,-3061)(0.00000,-60){11}{\makebox(1.6667,11.6667){\SetFigFont{5}{6}{\rmdefault}{\mddefault}{\updefault}.}}
}%
{\put(4201,-3661){\line(-1,-1){300}}
\put(3901,-3961){\line( 1,-1){300}}
\put(4201,-4261){\line(-1,-1){300}}
\put(3901,-4561){\line( 1,-1){300}}
\put(4201,-4861){\line( 1, 1){300}}
\put(4501,-4561){\line(-1, 1){300}}
\put(4201,-4261){\line( 1, 1){300}}
\put(4501,-3961){\line(-1, 1){300}}
}%
{\put(4501,-4561){\line( 1, 1){300}}
\put(4801,-4261){\line(-1, 1){300}}
\put(4501,-3961){\line( 1, 1){300}}
\put(4801,-3661){\line( 1,-1){300}}
\put(5101,-3961){\line(-1,-1){300}}
}%
{\multiput(5701,-3361)(-42.8571,-42.8571){15}{\makebox(1.6667,11.6667){\SetFigFont{5}{6}{\rmdefault}{\mddefault}{\updefault}.}}
}%
{\multiput(5701,-3361)(-42.8571,42.8571){15}{\makebox(1.6667,11.6667){\SetFigFont{5}{6}{\rmdefault}{\mddefault}{\updefault}.}}
}%
{\put(4201,-2461){\line( 1, 1){300}}
}%
{\put(4201,-3061){\line( 1, 1){600}}
}%
{\put(4201,-2461){\line( 1,-1){300}}
}%
{\put(4501,-2761){\line( 1,-1){300}}
}%
{\put(4801,-3061){\line( 1, 1){300}}
}%
{\put(4201,-1861){\line( 1,-1){900}}
}%
\thicklines
{\put(3901,-4561){\circle*{75}}
}%
{\put(4501,-4561){\circle*{75}}
}%
{\put(4801,-4261){\circle*{75}}
}%
{\put(4201,-4861){\circle*{75}}
}%
{\put(4201,-4261){\circle*{75}}
}%
{\put(3901,-3961){\circle*{75}}
}%
{\put(4201,-3661){\circle*{75}}
}%
{\put(4501,-3961){\circle*{75}}
}%
{\put(4801,-3661){\circle*{75}}
}%
{\put(5101,-3961){\circle*{75}}
}%
{\put(5701,-3361){\circle*{75}}
}%
{\put(5101,-2761){\circle*{75}}
}%
{\put(4801,-3061){\circle*{75}}
}%
{\put(4501,-2761){\circle*{75}}
}%
{\put(4801,-2461){\circle*{75}}
}%
{\put(4201,-3061){\circle*{75}}
}%
{\put(3901,-2761){\circle*{75}}
}%
{\put(3901,-2161){\circle*{75}}
}%
{\put(4201,-2461){\circle*{75}}
}%
{\put(4201,-1861){\circle*{75}}
}%
{\put(4501,-2161){\circle*{75}}
}%
\put(5801,-3381){$y$}%
\end{picture}$$

$$\text{Figure 2}$$

\begin{proof}
Observe that by Lemma 4, 
the elements of $g(C_k\ast S)$ lying strictly over $x_0$ are isomorphic to a
quotient of $g(C_{k-1} \ast S)$.  Now applying Lemma 3 inductively, we see
that every element of $g(C_k \ast S)$ can be written as $(y\vee x_i) 
\wedge x_j$ for $j\geq i$.  It follows that $g(C_k\ast S)$  is a quotient
of the lattice from Figure 2, but since the lattice from Figure 2 is graded,
it must coincide with $g(L)$. 
\end{proof}

\section*{Left Modular Lattices }

In this section, we recall a few results about left modular 
elements and left modular lattices from
[Li] and [MT].  

\begin{lemma}[{[Li]}] Suppose $u \lessdot v$ are left modular in $L$.
Let $z\in L$. Then:  

(i) $u\vee z \lessdoteq v \vee z$.

(ii) $u\wedge z \lessdoteq v \wedge z$.
\end{lemma}

\begin{proof} We prove (i).  Suppose otherwise, so that there is some element
$y$ such that $u\vee z < y < v\vee z$.  Now observe that
$((u\vee z)\vee v)\wedge y=y$. Now $v \wedge y  = u$, so 
$(u\vee z) \vee (v \wedge y)=u\vee z$, contradicting the left modularity
of $u$.  This proves (i).  Now (ii) follows by duality.  \end{proof}

\begin{lemma}[{[MT]}] Let $x$ be left modular, and $y<z$.  Then
$y\vee x \wedge z$ is left modular in $[y,z]$. \end{lemma}

\begin{proof} Let $s<t$ in $[y,z]$.  
$$(s\vee (y\vee x\wedge z))\wedge t= (s\vee x\wedge z)\wedge t=
s\vee x \wedge t= s\vee (y\vee x \wedge t)=s\vee((y\vee x\wedge z)\wedge t).$$
\end{proof}

\begin{lemma}[{[MT]}] If $L$ is a finite lattice with a maximal left 
modular chain $\zeroh=x_0\lessdot x_1\lessdot \dots \lessdot x_r=\oneh$, 
and $y\leq z$, then the set of elements of the form $y\vee x_i\wedge z$ forms a 
maximal left modular chain in $[y,z]$.  \end{lemma}

\begin{proof} 
The fact that the elements of the form 
$y\vee x_i\wedge z$ form a maximal chain in $[y,z]$ 
follows from
Lemma 5; the fact that they are left modular, from Lemma 6.
\end{proof}

\section*{Modularity}

For $y\leq z$, let us write $M(x,y,z)$ for the statement:
$$M(x,y,z):\  (y\vee x)\wedge z=y \vee (x \wedge z).$$
A lattice is said to be {\it modular} if $M(x,y,z)$ holds for all 
$x$ whenever
$y\leq z$.  

Standard notation is to write $xMz$ for the statement that $M(x,y,z)$ 
holds for all
$y\leq z$.  In this case $(x,z)$ is called a {\it modular pair}. 
An element $x$ is said to be 
modular if for any $z$ both $xMz$ and $zMx$ are modular pairs.
As we have already seen, an element $x$ is {\it left modular} if it
satisfies half the condition   
of being modular, namely that $xMz$ for all $z$.  

Let $L$ be a finite graded left modular lattice, with maximal 
left modular chain $\zeroh=x_0\lessdot x_1\lessdot\dots\lessdot x_r=\oneh$,
which we denote $\bx$.    
By definition, for any $y\leq z$, we have $M(x_i,y,z)$.  
We also have the following lemma:

\begin{lemma} In a finite graded left modular lattice $L$, with
maximal left modular chain $\bx$, for any $w\in L$ and $i<j$, we have
$M(w,x_i,x_j)$.  \end{lemma}

\begin{proof} 
Consider the sublattice $K$ of $L$ generated by $\bx$ and $w$.
First, we show that $K$ is graded.    
Let $y< z \in K$.  
By Lemma 7, we know that the elements of the form
$y\vee x_i \wedge z$ form a maximal chain in $L$.
These are all elements of $K$, so there is a maximal chain between
$y$ and $z$ having the same length as in $L$.  It follows that the 
covering relations in $K$ are a subset of the covering relations in $L$,
and hence that $K$ is graded (with the same rank function as $L$).  

Since $K$ is generated by $\bx$ and $w$, $K$ is a quotient of 
$C_r \ast S$.  Further, since $K$
is graded, it is a quotient of $g(C_r\ast S)$.  Since $g(C_r\ast S)$
is distributive, the modular equality is always satisfied in it, and
therefore also in $K$.  So $M(w,x_i,x_j)$ holds in $K$, and therefore
in $L$.  \end{proof}

\section*{Graded Left Modular Lattices are Supersolvable}

In this section, we prove Theorem 1, that finite graded 
left modular lattices are supersolvable. To do this, we have to show
that the sublattice generated by the left modular chain and another chain
is distributive.  

The proof mimics the proof of Proposition 2.1 of [St], which shows that if
$L$ is a finite lattice with a maximal 
chain of modular elements, 
then this chain is an $M$-chain, and hence $L$ is
supersolvable.  
The proof from [St] is based on Birkhoff's proof [Bi, \S III.7] 
that a modular lattice
generated by two chains is distributive.  

We recall briefly the way Birkhoff's proof works.  
Let $L$ be a finite modular 
lattice, and let
$\zeroh=x_0<\dots<x_r=\oneh$ and 
$\zeroh=y_0<\dots<y_s=\oneh$ be
two chains, which we denote $\bx$ and $\by$ respectively.  Assume further
that $L$ is generated by $\bx$ and $\by$.  
Let $u^i_j=x_i\wedge y_j$, and let
$v^i_j=x_i\vee y_j$.  Write $U$ for $\{u^i_j\}$ and $V$ for $\{v^i_j\}$.

Observe [Bi, \S III.7 Lemma 1] that any join of elements of $U$ can
be written in the form
$$\bigvee_{i=1}^t a_i\wedge b_i$$
where $a_1,a_2,\dots$ form a decreasing sequence from $\bx$, and 
$b_1,b_2,\dots$ form an
increasing sequence from $\by$.  

Then [Bi, \S III.7 Lemma 2], the following two 
identities are established for all decreasing sequences 
$a_1, a_2, \dots$ from $\bx$ and increasing sequences $b_1,b_2,\dots$ from
$\by$, 
and for all $t$, 
under the assumption
that $L$ is modular:

\begin{eqnarray*}
&{\bold P}_t\!:&({b_1}\vee {a_1})\wedge ({b_2}\vee {a_2}) \wedge \dots
\wedge ({b_t}\vee {a_t})=
{b_1}\vee ({a_1}\wedge {b_2}) \vee \dots
\vee ({a_{t-1}}\wedge {b_t})\vee {a_t} \\
&{\bold Q}_t\!:&({a_1}\wedge {b_1})\vee ({a_2}\wedge {b_2}) \vee \dots
\vee ({a_t}\wedge {b_t})=
{a_1}\wedge ({b_1}\vee {a_2}) \wedge \dots
\wedge ({b_{t-1}}\vee {a_t})\wedge {b_t} 
\end{eqnarray*}

Using $\bP_t$ and $\bQ_t$, it is straightforward 
to see that the set of joins of elements of $U$ 
coincides with the set of meets of elements of $V$
and that they therefore 
form a sublattice of $L$ [Bi, III\S7 Lemma 3].  Since $L$
is generated by $\bx$ and $\by$ (by assumption), it follows that
every element of $L$ can be written as a join of elements of $U$.

We now deviate slightly from the exposition in [Bi].  
Let $D$ denote the distributive lattice of down-closed subsets of 
$[1,r] \times [1,s]$.  Define a map $\phi:D \rightarrow L$ by
setting 
$$\phi(I)=\bigvee_{(i,j)\in I} u^i_j.$$
This map respects join operations, and from what we have already shown, 
it is surjective.  

Similarly, define a map $\psi:D\rightarrow L$ by setting
$$\psi(I)=\bigwedge_{(i,j)\not\in I} v^{i-1}_{j-1}.$$
This map respects meet operations.  Now, we observe
(by $\bP_t$ and $\bQ_t$) that $\phi$ and $\psi$ coincide.  They therefore
form a lattice homomorphism from $D$ onto $L$, which shows that $L$ is 
distributive, as desired.  

The only point at which modularity has been used is in 
establishing $\bP_t$ and $\bQ_t$.  Stanley noticed that it was sufficient
to assume only that all the $x_i$ are modular.  In fact, still less is 
sufficient.  

\begin{lemma} $\bP_t$ and $\bQ_t$ hold in any graded lattice such that
the $x_i$ form a maximal chain of left modular elements.  \end{lemma}

\begin{proof}
We prove $\bP_t$ and $\bQ_t$ by simultaneous induction on $t$.  
$\bP_1$ and $\bQ_1$ are tautologous.  Assume that $\bP_{t-1}$ and 
$\bQ_{t-1}$ hold.  
We now prove $\bQ_t$.  Recall that $a_1,a_2,\dots$ is a decreasing
sequence from $\bx$, and $b_1,b_2,\dots$ is an increasing sequence from $\by$.
We start from the lefthand side of $\bQ_t$:
\begin{eqnarray*}
& & ({a_1} \wedge {b_1})\vee \dots \vee({a_{t-1}}\wedge {b_{t-1}}) \vee 
 ({a_t}\wedge {b_t})\\
& &\qquad=\big(({a_1}\wedge {b_1})\vee \dots \vee ({a_{t-1}}\wedge 
{b_{t-1}}) \vee  {a_t}\big)\wedge {b_t}\big)\\
& & \qquad\qquad\text{by $M({a_t}, ({a_1}\wedge {b_1})\vee \dots \vee
({a_{t-1}}\wedge {b_{t-1}}),b_t)$ }\\
& & \qquad=\big[\big({a_1}\wedge({b_1}\vee {a_2})\wedge \dots \wedge (b_{t-2}\vee
a_{t-1})\wedge b_{t-1}\big) \vee a_t \big] \wedge b_t \\
& & \qquad\qquad\text{by $\bQ_{t-1}$}\\
& & \qquad=a_1 \wedge\big[\big(({b_1}\vee {a_2})\wedge \dots \wedge (b_{t-2}\vee
a_{t-1})\wedge b_{t-1}\big) \vee a_{t}\big] \wedge b_t \\
& & \qquad\qquad\text{by $M(({b_1}\vee {a_2})\wedge \dots \wedge (b_{t-2}\vee
a_{t-1})\wedge b_{t-1},a_t,a_1)$ (Lemma 8)}\\
& & \qquad=a_1 \wedge\big[\big(({b_1}\vee {a_2})\wedge \dots \wedge (b_{t-2}\vee
a_{t-1})\wedge (b_{t-1}\vee \zeroh)\big) \vee a_{t}\big] \wedge b_t \\
& & \qquad=a_1 \wedge\big[\big({b_1}\vee ({a_2}\wedge b_2) \dots 
\vee
(a_{t-1}\wedge b_{t-1})\big) \vee a_{t}\big] \wedge b_t \\
& & \qquad\qquad\text{by $\bP_{t-1}$}\\
& & \qquad=a_1\wedge\left[(b_1\vee a_2)\wedge\dots\wedge (b_{t-1}\vee a_t)\right] \wedge
b_t \\
& & \qquad\qquad\text{by $\bP_{t-1}$.}
\end{eqnarray*}

This proves $\bQ_t$.  
The dual argument holds for $\bP_t$, which completes the 
induction step, and the proof of the lemma
\end{proof}

This shows that Birkhoff's proof can be adapted to our situation, proving
Theorem 1.  

\section*{ Acknowledgements}

I would like to thank George Gr\"atzer for suggesting the idea for this 
paper, and for some guidance about free lattices.  I would like to thank
Peter McNamara, Ralph Freese, and an anonymous referee 
for helpful suggestions.


\begin{thebibliography}{MT}




\bibitem[Bi]{Bi} G. Birkhoff, \emph{Lattice Theory}, third edition.
American 
Mathematical Society, Providence, RI, 1973.

\bibitem[BS]{BS} A. Blass and B. Sagan, \emph{M\"obius functions of lattices}.
Adv. Math. \textbf{127} (1997), 94--123.

\bibitem[Gr]{Gr} G. Gr\"atzer, \emph{
General Lattice Theory}, 
second edition.
Birk\"auser, Basel, 2003.

\bibitem[Li]{Li} L. S.-C. Liu, \emph{Left-modular elements and edge labellings}.
Ph. D. thesis, Michigan State University, 1999.


\bibitem[Mc]{Mc} P. McNamara, \emph{EL-labelings, supersolvability, and 
0-Hecke algebra actions on posets}. J. Combin. Theory Ser. A \textbf{101} 
(2003), 69--89.

\bibitem[MT]{MT} P. McNamara and H. Thomas, \emph{Poset Edge-Labelling
and Left Modularity}. European J. Combin., to appear.  

\bibitem[St]{St} R. Stanley, \emph{Supersolvable lattices}.  
Algebra Universalis \textbf{2} (1972), 197--217.

 
\end{thebibliography}
\end{document}